\documentclass{svjour2}
\smartqed
\usepackage{epsfig}
\usepackage{graphicx}
\usepackage{amsmath,amsfonts,amssymb}

\begin{document}

\title{A Haar-like Construction for the Ornstein Uhlenbeck Process}

\author{Thibaud Taillefumier \and Marcelo O. Magnasco}

\institute{Thibaud Taillefumier \at
              Laboratory of Mathematical Physics, The Rockefeller University,10021 New York, NY USA\\
              \email{ttaillefum@rockefeller.edu}           
           \and
           Marcelo O. Magnasco \at
             Laboratory of Mathematical Physics, The Rockefeller University,10021 New York, NY USA\\
             \email{magnasco@rockefeller.edu}
}

\date{Received: date / Accepted: date}

\maketitle

\begin{abstract}
The classical Haar construction of Brownian motion uses a binary tree of triangular wedge-shaped functions. This basis has compactness properties which make it especially suited for certain classes of numerical algorithms. We present a similar basis for the Ornstein-Uhlenbeck process, in which the basis elements approach asymptotically the Haar functions as the index increases, and preserve the following properties of the Haar basis: all basis elements have compact support on an open interval with dyadic rational endpoints; these intervals are nested and become smaller for larger indices of the basis element, and for any dyadic rational, only a finite number of basis elements is nonzero at that number. Thus the expansion in our basis, when evaluated at a dyadic rational, terminates in a finite number of steps. We prove the covariance formulae for our expansion and discuss its statistical interpretation and connections to asymptotic scale invariance. 
\end{abstract}

\section{Introduction}

Random walks and continuous stochastic processes are of fundamental importance in a number of applied areas, including optics ÷\cite{Risken}, chemical physics ÷\cite{Kampen}, biophysics ÷\cite{Gerstner,Koch}, biology ÷\cite{Berg} and finance ÷\cite{Rolski}.
The mathematical idealization of the one-dimensionnal continuous random walk, the Wiener process, can be expressed in ---infinitely--- many bases as a sum of random coefficients times basis elements. Unique among these bases, the Haar ---or Schauder--- basis has three properties that make it particularly suitable for certain numerical computations. 
First, the basis elements all have {\em compact support}: the basis elements are nonzero only in open intervals. 
Second, the support is {\em increasingly compact}, i.e., the open intervals become smaller for higher indices of the basis elements; in fact, the intervals are nested in binary tree fashion, and have dyadic rational endpoints. 
Finally, given any dyadic rational, there is a finite number of basis elements which are nonzero at that number, so that evaluation of the Haar expansion at a dyadic rational {\em terminates} in a finite number of steps known beforehand. 
These properties can be used to great advantage in algorithms that construct the random walk in a ``top-down'' fashion, such as dychotomic search algorithms for first passage times. \\
However, the ``plain'' Wiener process has limited applicability in the areas mentioned above, so an extension of this construction to more complex stochastic processes is desirable. 
The naive generalization of the Haar basis construction to other stochastic processes fails to display our three properties. We present a method for constructing a Haar-like basis for the Ornstein-Uhlenbeck process which preserves these properties. The basis is therefore useful for advanced numerical computations: a fast dichotomic search algorithm for first passage time computations shall be presented elsewhere. The method we present is also amenable to further generalizations to other stochastic processes. \\
This paper is organized as follows. We first review some background in stochastic processes and basis expansions. Then we review the well-known decomposition of a Wiener process in a basis of functions derived from the Haar system. 
In the third section, we give a statistical interpretation of such a construction, which leads us to propose a basis for the Ornstein-Uhlenbeck process.
In the fourth section, we  prove that the Ornstein-Uhlenbeck process is correctly represented as a discrete process in the proposed basis.
In the last section, we further the statistical interpretation and its connection to scale invariance and Markovian properties.

\section{Background on Stochastic Processes}

The Wiener process and the Ornstein-Uhlenbeck process are continuous stochastic processes; 
we specify this class of process through the Langevin equation
\begin{equation} \label{eq:stochastic system}
\left\{ \begin{array}{l}
\dot{x} = f(x,t)+ \eta(t) \\
x(t_{0}) = x_{0}, \: t \in \left[ t_{0},T\right] \: ,
\end{array} \right.
\end{equation}
where $f$ is a deterministic function and $\eta(t)$ describes the stochastic forcing.
Equation \eqref{eq:stochastic system} is a first order stochastic differential equation and its connection to the Fokker-Planck equation has been extensively studied ÷\cite{Risken,Kampen}. 
We only consider here the case where the noise is white and Gaussian: $\eta(t)$ are realizations of independent identically distributed Gaussian variables $\eta_t$, with time correlations satisfying
\begin{equation}
\langle \, \eta_{t} \cdot \eta_{s} \, \rangle = \Gamma \cdot \delta (s-t) \: ,
\nonumber
\end{equation}
where $\delta$ is the Dirac distribution.\\
We denote by $\omega$ a given realization of the stochastic forcing: the collection of all the values  ${\lbrace \eta(t) \rbrace }_{t \in \left[ t_{0},T \right]}$ in an interval. 
The set of $\omega$ values defines the sample space $\Omega$ and the probability of occurence of a sequence $\omega$ in $\Omega$ is
determined by the joint probability of ${\lbrace \eta_{t} \rbrace}_{t \in \left[ t_{0},T\right]}$. 
With this notation, we can introduce the general solution of the stochastic system as the stochastic process $X_{t}$.
For a given realization of the noise $\omega$, there is a unique solution to \eqref{eq:stochastic system} called a sample path: neglecting to notate the dependence on the initial condition, 
we write $X_{t}(\omega)$ the value of this sample path at time $t$.
Since each sample path ${\lbrace X_{t}(\omega) \rbrace}_{t \in \left[ t_{0},T\right]}$ occurs with the same probability as its matching sequence $\omega$ in the sample space $\Omega$, the value $X_{t}(\omega)$ can be seen as the outcome of a random function $X_{t}$  defined on $\Omega$. 
$X_{t}$ is the stochastic process solution of  \eqref{eq:stochastic system} and has several important properties: it is a continuous process as it is defined for a continuous index set  $\left[ t_{0},T\right]$; it is a Gaussian process, as it integrates contributions of Gaussian variables; and, being a Markovian process, the value of $X_t$ only depends on ${\lbrace \eta(u) \rbrace }_{u \in \left[ 0,t \right]}$, the sequence of realizations preceding $t$. 
Two special forms of $f$ shall concern us. When $f$ is zero, the process is called the Wiener process $W_{t}$; when $f$ is linear in $x$, the process is the Ornstein-Uhlenbeck process $U_{t}$. Due to the relative simplicity of both situations, the probability laws of the processes 
(i.e., the Green functions of the associated Fokker-Planck equations) are known analytically. 
If a Wiener process is at $x_{0}$ at time $t=t_{0}$, the probability of finding the process in $x$ at time $t$ is
\begin{equation} \label{eq:probRW}
\mathrm{ \bold P}( W_{t} \scriptstyle = \displaystyle \! x \, \vert W_{t_{0}} \scriptstyle = \displaystyle \! x_{0} ) = \frac{1}{\sqrt{2 \pi} \cdot {}_{ \scriptscriptstyle W \displaystyle} \sigma_{t}} \cdot \exp{ \bigg( - \! \frac{ \left( x-x_{0} \right)^{2}}{2 \cdot {}_{ \scriptscriptstyle W \displaystyle} \sigma_{t}^{2}} \bigg) } 
\end{equation}
with a variance ${}_{ \scriptscriptstyle W \displaystyle} \sigma_{t}^{2} = \Gamma \! \cdot \! \left( t - t_0 \right)$. 
For the Ornstein-Uhlenbeck process, a similar result holds 
\begin{eqnarray} \label{eq:probOU}
\mathrm{ \bold P}( U_{t} \scriptstyle = \displaystyle \! x \, \vert U_{t_0} \scriptstyle = \displaystyle \! x_{0} )  = 
 \frac{1}{\sqrt{2 \pi} \cdot {}_{ \scriptscriptstyle U \displaystyle} \sigma_{t}} \cdot \exp{ \left( -  \frac{ \left( x \! - \! x_{0} \, e^{-\alpha \left( t -  t_0 \right)} \right)^{2}}{2 \cdot {}_{ \scriptscriptstyle U \displaystyle} \sigma_{t}^{2}} \right) }
\end{eqnarray}
with a variance ${}_{ \scriptscriptstyle U \displaystyle} \sigma_{t}^{2} = \frac{\Gamma}{2 \alpha} \! \cdot \! \left( 1 - e^{-\alpha \left( t - t_0 \right)} \right)$. 
The previous expressions describe the statistics of $W_t$ and $U_t$, which will be called $X_t$ when collectively designated. 
Continuous processes require an infinite number of random variables, and establishing results about them is quite difficult.
Consider for example the Ornstein-Uhlenbeck process,  widely used from finance to neuroscience: finding analytically the first-passage times distribution with a fixed threshold proves a surprisingly intricate question in this situation ÷\cite{Siegert,Ricciardi,Leblanc}, and numerically, sample paths are only approximated by stochastic Euler methods, with integration schema of low efficiency ÷\cite{Giraudo,Bouleau,Gardiner}.
Abating these difficulties for the Ornstein-Uhlenbeck process is the motivation for this paper.\\
To circumvent the problem, it is advantageous to represent a continuous process as a discrete process. 
Conspicuously enough, a discrete process has a countable index set of random variables. 
At stake is to write a Gaussian process $X_{t}$ as a convergent series of random functions $f_{n} \cdot \xi_{n}$, where $f_{n}$ is a deterministic function and $\xi_{n}$ a Gaussian variable of law $\mathcal{N}(0,1)$ ( i.e.\!\! with null mean and unitary variance ). 
Assuming the coefficients of the decomposition to be included in the definition of $f_{n}$, the identity 
\begin{equation}
X_{t} = \sum_{n=0}^{\infty} f_{n} (t)\cdot \xi_{n}  = \lim_{N \to \infty} \sum_{n=0}^{N} f_{n} (t)\cdot \xi_{n} 
\nonumber
\end{equation}
shows $X_{t}$ as the limit of a sequence of finite processes $\sum_{n=0}^{N} f_{n} (t)\cdot \xi_{n}$. 
Depending on the nature of the convergence, this may result in two advantages. 
Analytically, it is generally more tractable to prove mathematical results on finite combination of simple random functions and then to extend them to a limit random process. 
Numerically, the quantity $\sum_{n=0}^{N} f_{n} (t)\cdot \xi_{n}$ can be accurately computed and gives a correct approximation of the process at any level of precision. \\


\section{The Haar construction of the Wiener process} 

\begin{figure}
\begin{center}
\vspace{-40pt}
\includegraphics[scale=0.55]{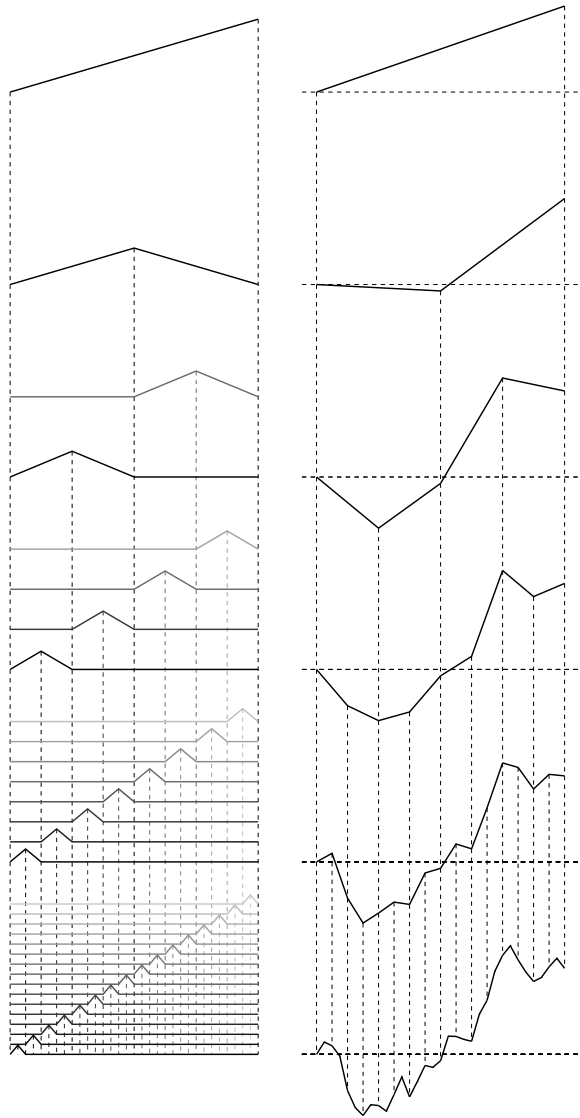}
\caption{In the left column, the elements of the basis $\Psi_{n,k}$ are represented for each rank $n$ with $0 \! \leq \! n \! < \! 6$. In the right column, the partial sums $W^{n}(\omega)$ are shown for a given set of realizations $\omega$. Note that each element $\Psi_{n,k}$ has a compact support delimited by dyadic numbers in $D_{n} = \lbrace k2^{-n} \, \vert \, 0  \! \leq \! k \! \leq \! 2^{n}\rbrace$ and that all $\Psi_{n',k}$  is zero on $D_{n}$ for $n'>n$.}
\label{fig:RWbasis}
\end{center}
\end{figure}

The Haar system is the set of functions $h_{n,k}$ in $L^{2}([0,1])$ defined by 
\begin{equation} \label{eq:haar definition}
h_{n,k}(t)= \left\{ \begin{array}{lll}
\! & \! \! 2^{\frac{n-1}{2}} & \textrm{if $\left(2k\right)2^{-n} \! \leq \! t \! \leq \! \left(2k\!+\!1\right)2^{-n} $} \: ,\\
\! -& \! \! 2^{\frac{n-1}{2}} & \textrm{if $\left(2k\!+\!1\right)2^{-n} \! \leq \! t \! \leq \! 2\left(k\!+\!1\right)2^{-n}$} \: ,\\
\! & \! \! 0  & \textrm{otherwise} \: . \end{array} \right. 
\nonumber\\
\end{equation}
for $n\!\geq\!1$ with the addition of the function $h_{0,0}(t)\!=\!1$ on $[0,1]$. 
The Haar system has a several interesting properties.
First, the functions $h_{n,k}$ form a complete orthonormal basis of $L^{2}([0,1])$ for the scalar product $\left( f , g \right) = \int_{0}^{1} f(t)g(t)dt$. 
Second, each element $h_{n,k}$ has a compact support 
\begin{equation}
S_{n,k}\!= \! \left[ k \! \cdot \!2^{-n\!+\!1}, (k\!+\!1)2^{-n\!+\!1}\right]
\nonumber
\end{equation}
and, for a given $n$, the collection of supports $S_{n,k}$ represents a partition of $\left[ 0,1\right]$. 
Third, the functions $h_{n,k}$ build up a wavelet basis of $L^{2}([0,1])$, since we have the scale-invariant construction rule 
\begin{equation}
h_{n,k}(t) = 2^{\frac{n-1}{2}}  \cdot h_{1,0}(2^{n-1} t - k) \: .
\end{equation} 
Such properties prove useful to decompose simple Gaussian process as related in the following.
Consider the Wiener process $W_t$ as the stochastic integral of the independent random variables $\eta_{t}$ which follow a normal law $\mathcal{N}(0,\sqrt{\Gamma})$.
We introduce the associated Gaussian white noise process formally defined as $\frac{dW}{dt}$ on $\Omega$.
A sample path $\frac{dW_{t}}{dt} (\omega)$ is almost surely an element of $L^{2}([0,1])$ for a given noise realization  $\omega = {\lbrace \eta(t) \rbrace }_{t \in \left[ 0,1 \right]}$.
We write its decomposition on the Haar basis 
\begin{equation}
\left( \frac{dW_{t}}{dt} \right) (\omega) = \sum_{n=0}^{\infty} \sum_{ 0\leq k<2^{n\!-\!1} } c_{n,k}(\omega)h_{n,k}(t) \: ,
\nonumber
\end{equation}
introducing  $c_{n,k}(\omega)$ the component of $\frac{dW_{t}}{dt}(\omega)$ in the direction of $h_{n,k}$.
Each sequence $\omega$ can yield different coefficients $c_{n,k}(\omega)$ and their values appear as the outcome of a random variable $c_{n,k}$ defined on the sample space $\Omega$ 
\begin{equation}
c_{n,k} = \int_{0}^{1} h_{n,k}(t)\frac{dW_t}{dt} \, dt = \int_{0}^{1} h_{n,k}(t)dW_t \: .
\end{equation}
In our specific case of Gaussian uncorrelated noise, stochastic integration shows that the random variables  $c_{n,k}$ are all independent and identically distributed following the law $\mathcal{N}(0,\sqrt{\Gamma})$. 
The white noise process $\frac{dW_{t}}{dt}$ is then expressed as a discrete process on $\Omega$ by 
\begin{equation}
\frac{dW_{t}}{dt} = \sum_{n=0}^{\infty} \sum_{  0\leq k<2^{n\!-\!1} } \sqrt{\Gamma} \cdot h_{n,k}(t) \cdot \xi_{n,k} \: ,
\nonumber
\end{equation} 
where the $\xi_{n,k}$ are independent and distributed with law $\mathcal{N}(0,1)$.\\
The use of Haar functions to decompose white Gaussian noise directly suggests a corresponding result for the Wiener process.
The Wiener process $W_{t}$ is the stochastic integral of the Gaussian white noise process $\frac{dW_{t}}{dt}$, which leads to introduce the integrals of $h_{n,k}$ as candidates to build a basis of functions for the Wiener process ÷\cite{Karatzas,Pyke}.
We note this set of integral functions as
\begin{equation} \label{eq:Psi definition}
\Psi_{n,k}(t) = \sqrt{\Gamma} \cdot \int_{0}^{1} \chi_{ \left[ 0, t \right] }(u)h_{n,k}(u) \, du \: ,
\end{equation}
with the help of the indicator functions given by 
\begin{equation} \label{eq:indicative definition}
\chi_{ \left[ 0, t \right] }(u)= \left\{ \begin{array}{ll}
1& \textrm{if $0 \leq u \leq t $}\\
0  & \textrm{otherwise} \quad \end{array} \right. \: .
\nonumber
\end{equation}
The first elements of the so-defined basis are shown on figure \ref{fig:RWbasis}.
The question is then to know whether the process $W_{t}^{N}$ defined as the finite sum of random functions 
\begin{equation} \label{eq: Brownian sum}
W_{t}^{N}(\omega) = \sum_{n=0}^{N} \sum_{  0\leq k<2^{n\!-\!1} } \Psi_{n,k}(t) \cdot \xi_{n,k}(\omega)
\nonumber
\end{equation}
converges toward a Wiener process. 
It can be shown that the series converges normally to a limit process almost surely on $\Omega$  ÷\cite{Karatzas,Pyke}. 
Even though we have not proven its Wiener process nature yet,  we refer to this limit as $W_{t}$. 
Due to the normal convergence, $W_{t}$ is continuous in $t$ and being a sum of Gaussian variables, it is a Gaussian process.
Therefore, showing that $W_{t}$ is a Wiener process only amounts to demonstrate it has the same law of covariance as a Wiener process  ÷\cite{Karatzas,Pyke}, i.e. $ \langle W_{t} \cdot W_{s} \rangle =\Gamma \cdot \, \min\left( t , s \right) $, where $\min\left(t , s \right)$ is the minimum of $t$ and $s$. In other words, we need to evaluate the quantity
\begin{eqnarray} \label{eq: variance Brownian sum 0}
 \langle W_{t} \cdot W_{s} \rangle = \lim_{N \to \infty} \langle W_{t}^{N} \cdot W_{s}^{N} \rangle = 
 \lim_{N \to \infty}  \sum_{n=0}^{N} \sum_{ 0\leq k<2^{n\!-\!1} } \Psi_{n,k}(t) \cdot \Psi_{n,k}(s)
\nonumber
\end{eqnarray}
which entails the calculation of a rather tedious series. Indeed, for a given $t$, at each step $n$, there is only one $k$ for which $\Psi_{n,k}(t)$ is nonzero, and expressing the series analytically results in a complicated operation. One way to overcome the issue is to notice that the expected covariance result can be expressed in terms of
\begin{equation} \label{eq:wedge}
\min\left( t , s \right) = \int_{0}^{1} \chi_{ \left[ 0, t \right] }(u)\chi_{ \left[ 0, s \right] }(u) \, du   \: .
\end{equation}
The right term of \eqref{eq:wedge} is actually the scalar product of the functions $\chi_{ \left[ 0, t \right] }$ and $\chi_{ \left[ 0, s \right] }$ on $L^{2}([0,1])$. The key point is then to introduce
the decomposition in the Haar orthonormal basis to write the scalar product of two given functions as
\begin{eqnarray} \label{eq:haar decomposition}
 \int_{0}^{1} f(t)g(t) \, dt   \; =
 \sum_{\substack{n\geq0 \\  0\leq k<2^{n\!-\!1}}} \int_{0}^{1} f(u)h_{n,k}(u) \, du \int_{0}^{1} g(u)h_{n,k}(u) \, du \:  .
\end{eqnarray}
When applied to the indicator functions of interest $\chi_{ \left[ 0, t \right] }$ and $\chi_{ \left[ 0, s \right] }$, the relation \eqref{eq:haar decomposition} leads to
\begin{eqnarray} \label{eq:Psi decomposition}
\Gamma \cdot \min\left(t,s\right)  & = & \Gamma \cdot \int_{0}^{1} \chi_{ \left[ 0, t \right] }(u)\chi_{ \left[ 0, s \right] }(u) \, du   
\nonumber \\
& = &  \sum_{n=0}^{\infty} \sum_{ 0\leq k<2^{n\!-\!1}} \Psi_{n,k}(t) \Psi_{n,k}(s) \:,
\nonumber
\end{eqnarray}
since  the definition \eqref{eq:Psi definition} describes $\Psi_{n,k}(t)$ as the coefficient relative to $h_{n,k}$ in the decomposition of  $\chi_{ \left[ 0, t \right] }$ on the Haar system. 
We can finally recap the result
\begin{eqnarray} \label{eq: variance Brownian sum 1}
\langle W_{t} \cdot W_{s}\rangle & = & \lim_{N \to \infty} \langle W_{t}^{N} \cdot W_{s}^{N} \rangle 
\nonumber\\ 
& = &  \lim_{N \to \infty}  \sum_{n=0}^{N} \sum_{ 0\leq k<2^{n\!-\!1}} \Psi_{n,k}(t) \cdot \Psi_{n,k}(s) \quad =  \quad \Gamma \cdot \min\left(t,s\right) \: ,
\nonumber
\end{eqnarray}
establishing the discrete description of the Wiener process as a normally convergent series of terms $\Psi_{n,k} \cdot \xi_{n,k}$, where $\Psi_{n,k} $ is a Haar-derived function and $\xi_{n,k}$ a random variable of normal law $\mathcal{N} \left( 0,1\right)$.


\section{Comparison of the Wiener and Ornstein-Uhlenbeck processes}

\begin{figure}
\begin{center}
\includegraphics[scale=0.5]{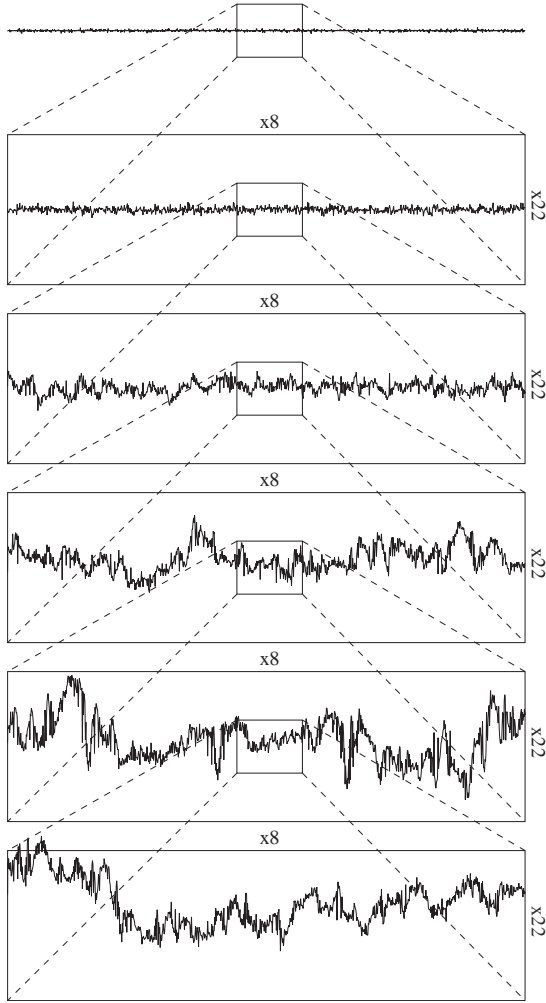}
\caption{A sample path ${U_{t}(\omega)}$ is represented at different magnifications following the scale invariance of a Wiener process: the vertical zooming factor is the square root of the horizontal factor. Note that the sample path ${U_{t}(\omega)}$ behaves as  a Wiener process at small scales.}
\label{fig:ScaleOU}
\end{center}
\end{figure}

We recall that the Langevin equation \eqref{eq:stochastic system} can be solved by quadratures in simple cases.
If the process is at $x_{0}$ when $t\!=\!0$, the Ornstein-Uhlenbeck process $U_{t}$ is expressed
\begin{equation} \label{eq:OUexpr}
U_{t} = x_{0} \, e^{-\alpha t} + \int_{0}^{t} e^{\alpha \left( u - t \right)} \eta(u) \, du \: ,
\end{equation}
as opposed to the Wiener process $W_{t}$ in the same conditions
\begin{equation} \label{eq:RWexpr}
W_{t} = x_{0} + \int_{0}^{t} \eta(u) \, du \: .
\end{equation}
The comparison of definitions \eqref{eq:OUexpr} and \eqref{eq:RWexpr} explains why finding a basis of decomposition for $U_{t}$ stumbles on a several difficulties.
First the process $U_{t}$ is not anymore a simple integral of white Gaussian noise, which is naturally described as a discrete process. 
Second, there is no more scale invariance, implying that a putative basis of decomposition is not to be thought of as wavelets. 
Finally, the exponential factor in \eqref{eq:OUexpr} indicates that the process $U_{t}$ does not sum the $\eta_{s}$ evenly; their contribution depends on the position of $s$ compared to $t$. 
As a consequence, the presence of these correlations makes it unlikely for the basis to conserve any orthogonality properties.
 \\
Yet, as noticeable in figure \ref{fig:ScaleOU}, the examination of a sample path ${\lbrace U_{t}(\omega) \rbrace}_{t \in \lbrack 0,1\rbrack}$ reveals the scale-invariant behavior of a Wiener process for asymptotically small time scale as well as for asymptotically small $\alpha$. 
It means that the basis of decomposition $\Psi_{n,k}$ for the Wiener process is asymptotically valid to describe $U_{t}$ at fine scale. 
This observation suggests that, upon slight alteration of its analytical expression, the Haar derived basis $\Psi_{n,k}$ can give rise to a basis $\Phi_{n,k}$ adapted to the Ornstein-Uhlenbeck process. 
The change in the analytical expression of $\Psi_{n,k}$ should be consistent with the previously mentioned difficulties, preventing its formulation to be scale invariant or orthogonal. 
Under this restraint, the fundamental property that each element $\Psi_{n,k}$ exhibits a compact support of the form $S_{n,k}$ should be preserved in the expression of $\Phi_{n,k}$.
\\
To carry out this program, the key point is to consider $\mathrm{ \bold P}( X_{t_{y}} \! \scriptstyle = \displaystyle \! y \, \vert X_{t_{x}} \! \scriptstyle = \displaystyle \! x, X_{t_{z}} \! \scriptstyle = \displaystyle \! z)$ with $t_{x}\!<\!t_{y}\!<\!t_{z}$, the probability law of $X_{t}$ knowing its values $x$ and  $z$ at two framing times $t_{x}$ and $t_{z}$. 
Because $X_t$ is a markovian process, a sample path ${\lbrace X_{t}(\omega) \rbrace}_{t \in \lbrack 0,1\rbrack}$ which originates from $x$ and joins $z$ through $y$ is just the junction of two independent paths: a path originating in $x$ going to $y$ and a path originating from $y$ going to $z$. 
Assuming conditional knowledge of its origin $x$, the probability of such a compound path is the product of the probability of the two elementary paths with conditional knowledge of their respective origins $x$ and $y$.
Therefore, after normalization by the absolute probability for a path to go from $x$ to $y$, $\mathrm{ \bold P}( X_{t_{y}} \! \scriptstyle = \displaystyle \! y \, \vert X_{t_{x}} \! \scriptstyle = \displaystyle \! x, X_{t_{z}} \! \scriptstyle = \displaystyle \! z)$ is expressed in the following expression
\begin{eqnarray} \label{eq:probrelation}
\mathrm{ \bold P}( X_{t_{y}} \! \scriptstyle = \displaystyle \! y \, \vert X_{t_{x}} \! \scriptstyle = \displaystyle \! x, X_{t_{z}} \! \scriptstyle = \displaystyle \! z) \quad 
= \quad
\frac{ \mathrm{\bold P}( X_{t_{y}} \! \scriptstyle = \displaystyle \! y\,  \vert X_{t_{x}} \! \scriptstyle = \displaystyle \! x) \cdot \mathrm{\bold P}( X_{t_{z}} \! \scriptstyle = \displaystyle \! z \, \vert X_{t_{y}} \! \scriptstyle = \displaystyle \! y) }{ \mathrm{\bold P}( X_{t_{z}} \! \scriptstyle = \displaystyle \! z \, \vert X_{t_{x}} \! \scriptstyle = \displaystyle \! x ) } \: .
\end{eqnarray}
It is now a simple matter of calculation to compute the distribution of  $X_{t_{y}}$ knowing $X_{t_x}\!\!=\!x$ and $X_{t_z}\!\!=\!z$ with the analytical expression of the probability $ \mathrm{\bold P}( X_{t_{y}} \! \scriptstyle = \displaystyle \! y\,  \vert X_{t_{x}} \! \scriptstyle = \displaystyle \! x)$.
In the case of a Gaussian process, it is expected to follow a normal law, which we refer to as $\mathcal{N}({}_{ \scriptscriptstyle P \displaystyle}\mu(t_{y}),{}_{ \scriptscriptstyle P \displaystyle}\sigma(t_{y}))$. 
For the Wiener process, using expression \eqref{eq:probRW} for $ \mathrm{\bold P}( X_{t_{y}} \! \scriptstyle = \displaystyle \! y\,  \vert X_{t_{x}} \! \scriptstyle = \displaystyle \! x)$, the mean value ${}_{ \scriptscriptstyle W \displaystyle} \mu(t_{y})$ and the variance ${}_{ \scriptscriptstyle W \displaystyle} \sigma(t_{y})^{2}$ result in
\begin{equation} \label{eq:muRW}
{}_{ \scriptscriptstyle W \displaystyle} \mu(t_{y}) = \frac{ t_{z}-t_{y} }{ t_{z}-t_{x} }  \cdot x + \frac{ t_{y}-t_{x} }{ t_{z}-t_{x} } \cdot z \: ,
\end{equation}
\begin{equation} \label{eq:sigmaRW}
{}_{ \scriptscriptstyle W \displaystyle} \sigma(t_{y})^{2} = \Gamma \cdot \frac{ (t_{y}-t_{x}) (t_{z}-t_{y}) }{ (t_{z}-t_{x}) } \: .
\end{equation}
For the Ornstein-Uhlenbeck process, using expression \eqref{eq:probOU} for $ \mathrm{\bold P}( X_{t_{y}} \! \scriptstyle = \displaystyle \! y\,  \vert X_{t_{x}} \! \scriptstyle = \displaystyle \! x)$ similarly yields the mean ${}_{ \scriptscriptstyle U \displaystyle} \mu(t_{y})$ and the variance ${}_{ \scriptscriptstyle U \displaystyle} \sigma(t_{y})^{2}$:
\begin{equation} \label{eq:muOU}
{}_{ \scriptscriptstyle U \displaystyle} \mu(t_{y}) = \frac{\sinh{ \big( \alpha(t_{z} \! - \! t_{y}) \big) }}{\sinh{\big( \alpha(t_{z} \! - \! t_{x}) \big) }}  \cdot x + \frac{\sinh{ \big( \alpha(t_{y} \! - \! t_{x}) \big) }}{\sinh{\big( \alpha(t_{z} \! - \! t_{x}) \big) }} \cdot z \: ,
\end{equation}
\begin{equation} \label{eq:sigmaOU}
{}_{ \scriptscriptstyle U \displaystyle} \sigma(t_{y})^{2} = \frac{\Gamma}{2 \alpha} \cdot  \frac{2 \cdot \sinh{ \big( \alpha(t_{y} \! - \! t_{x}) \big) } \cdot \sinh{ \big( \alpha(t_{z} \! - \! t_{y}) \big) }}{\sinh{\big( \alpha(t_{z} \! - \! t_{x}) \big) }} \: .
\end{equation}
In the limit of very short time scale or vanishing $\alpha$, we notice that  ${}_{ \scriptscriptstyle U \displaystyle} \mu(t_{y})$ and ${}_{ \scriptscriptstyle U \displaystyle} \sigma(t_{y})^{2}$ approximate ${}_{ \scriptscriptstyle W \displaystyle} \mu(t_{y})$ and  ${}_{ \scriptscriptstyle W \displaystyle} \sigma(t_{y})^{2}$.\\
We note $D_{N}$ the set of reals $\lbrace k2^{-N} \, \vert \, 0  \! \leq \! k \! \leq \! 2^{N}\rbrace$ and we have $\lbrace 0,1\rbrace = D_{0} \subset D_{1} \subset \cdots \subset D_{N}$  a growing sequence of sets with limit ensemble $\mathcal{D}$ the set of dyadic points in $\left[ 0, 1 \right]$.
Assuming we know the values of the process on a subset of dyadic points $D_{ \! N}$, we can construct the conditional average ${\langle P_{t} \rangle}_{D_{ \! N}}$, which is the most probable outcome of $X_{t}$ knowing its values on $D_{ N}$. For a Wiener process, \eqref{eq:muRW} shows that ${\langle W_{t} \rangle}_{D_{ \! N}}$ is a piece-wise linear function of t interpolating each points of $D_{ \! N}$; whereas for an Ornstein-Uhlenbeck process,  \eqref{eq:muOU} depicts  ${\langle U_{t} \rangle}_{D_{ \! N}}$ as a succession of catenaries joining successive points of $D_{ \! N}$.
With $0 \! \leq \! k \! < \! 2^{-N}$, if  $t_{x}\!=\!k2^{-N}$ and  $t_{z}\!=\!(k\!+\!1)2^{-N}$ are the two successive points of $D_{N}$ framing t, the average ${\langle X_{t} \rangle}_{D_{ \! N}}$ is only conditioned by $X_{t_{x}} \! \! = \! x$ and $X_{t_{z}} \! \! = \! z$. 
For the sake of simplicity, we write
\begin{equation}
{\langle X_{t} \rangle}_{D_{ \! N}} = {\langle X_{t} \rangle}_{x,z} = {}_{ \scriptscriptstyle X \displaystyle} \mu_{t_{x},t_{z}}(t, x, z) \stackrel{def}{=} {}_{ \scriptscriptstyle X \displaystyle} \mu^{N,k}(t)\: ,
\end{equation}
where the conditional dependency upon $x$ and $z$ is implicit in ${}_{ \scriptscriptstyle X \displaystyle} \mu^{N,k}$.
We want to investigate the change in the estimation of $P_{t}$ due to the conditional knowledge of its value on the dyadic set $D_{N+1}$. In that perspective, we exemplified the conditional average ${\langle X_{t} \rangle}_{D_{ \! N+\!1}}$ on $\left[ t_{x}, t_{z}\right]$ where the estimation of $X_{t}$ is now dependent upon the value $X_{t_{y}} \! \! = \! y$ with
$t_{y}$ the midpoint of $t_{x}$ and $t_{z}$:
\begin{eqnarray}\label{eq:muN+1}
{\langle X_{t} \rangle}_{D_{ \! N \! + \! 1}}   =  {\langle X_{t} \rangle}_{x,y,z} & = & 
 \left\{
\begin{array}{ll}
 {\langle X_{t} \rangle}_{x,y} \quad \textrm{if}  \; t_{x} \! \leq \! t \! \leq \! t_{y} \: ,\\
 {\langle X_{t} \rangle}_{y,z} \quad \textrm{if} \; t_{y} \! \leq \! t \! \leq \! t_{z} \: ,\\
\end{array} \right.
\nonumber\\
&  \stackrel{def}{=} & {}_{ \scriptscriptstyle X \displaystyle} \nu^{N,k}(t , y) \: .
\end{eqnarray}
We remark that, being a function of $y$, the conditional average ${\langle X_{t} \rangle}_{x,y,z}$ determines a random function ${}_{ \scriptscriptstyle X \displaystyle} \nu^{N,k}(t , Y_{N,k})$, where the short notation $Y_{N,k}$ indicates the Gaussian variable $X_{t_y}$ knowing $X_{t_x}\!\!=\!x$ and $X_{t_z}\!\!=\!z$. The probability distribution of $Y_{N,k}$ follows the law $\mathcal{N}({}_{ \scriptscriptstyle X \displaystyle}\mu(t_{y}),{}_{ \scriptscriptstyle X \displaystyle}\sigma(t_{y}))$ and it gives, through the function ${}_{ \scriptscriptstyle X \displaystyle} \nu^{N,k}$, the random contribution of ignoring $X_{t_y}\!\!=\!y$ when one estimates the process knowing its values on $t_x$ and $t_z$.
\\
The results above allows to gain insight in the building of a Wiener process $W_{t}$ as the converging series of random functions $\Psi_{n,k} \cdot \xi_{n,k}$.
It is easy to see from the definition \eqref{eq:haar definition} that $\Psi_{n,k}$ is linear between any two points in $D_{n}$ for $n \leq N $ and that  $\Psi_{n,k}$ is zero on $D_{n}$ for every $n > N $. In other words, the partial sum 
\begin{equation}
W_{t}^{N} = \sum_{n=0}^{N} \sum_{ 0\leq k<2^{n\!-\!1} } \Psi_{n,k}(t) \cdot \xi_{n,k}\quad \mathrm{for} \quad t \in D_{N}  
\nonumber
\end{equation}
coincide with $W_{t}$ on $D_{N}$ and more generally with ${\langle X_{t} \rangle}_{D_{ \! N}}$ on $\left[ 0, 1 \right]$. Identifying partial sums with conditional averages, it is then straightforward to express the component $\Psi_{N\!+\!1,k}(t) \cdot \xi_{N\!+\!1,k}$ in the decomposition of $W_{t}$ 
\begin{eqnarray}
\Psi_{N\!+\!1,k}(t) \cdot \xi_{N\!+\!1,k} & = & W_{t}^{N\!+\!1} - W_{t}^{N} \hspace{30pt} \nonumber\\
& = & {\langle W_{t} \rangle}_{D_{ \! N \! + \! 1}} - {\langle W_{t} \rangle}_{D_{ \! N}} \\
& = & {}_{ \scriptscriptstyle W \displaystyle} \nu^{N,k}(t ,Y_{N,k}) - {}_{ \scriptscriptstyle W \displaystyle} \mu^{N,k}(t) \: , \nonumber
\end{eqnarray} 
bearing in mind the previous definitions for which $\left[ t_{x}, t_{z}\right] \! = \! \left[ k2^{-N}, (k\!+\!1)2^{-N}\right]$ is the support $S_{N\!+\!1,k }$ of $\Psi_{N\!+\!1,k}$. The tight connection between $\xi_{N\!+\!1,k}$ and $Y_{N,k}$ is made obvious: if one knows the values of the process on $D_{N}$, the random contribution of $\sum_{k}\Psi_{N\!+\!1,k}(t) \cdot \xi_{N\!+\!1,k}$ conveys the uncertainty about $W_t$ that is discarded by the knowledge of its values on $D_{N\!+\!1} \setminus D_{\!N}$. \\
We are now in a position to complete our program: continuing the identification of partial sums and conditional average for the Ornstein-Uhlenbeck process $U_t$, it is direct to propose a basis of decomposition
\begin{equation} \label{eq:defPsi}
\Phi_{N\!+\!1,k}(t) \cdot \xi_{N\!+\!1,k}  =  {}_{ \scriptscriptstyle U \displaystyle} \nu^{N,k}(t ,Y_{N,k}) - {}_{ \scriptscriptstyle U \displaystyle} \mu^{N,k}(t) 
\end{equation} 
 with the adapted definitions of ${}_{ \scriptscriptstyle U \displaystyle} \mu^{N,k}$ and ${}_{ \scriptscriptstyle U \displaystyle} \nu^{N,k}$ on $S_{N\!+\!1,k}$, the support of the investigated functions $\Phi_{N\!+\!1,k}$. We underline that the notation $Y_{N,k}$ refers here to the random process $U_t$ at the midpoint of the support $t=(2k\!+\!1)2^{-(N+1)}$ knowing its values on the extremities.


\section{A discrete basis of functions to generate the Ornstein-Uhlenbeck process}

\begin{figure}
\begin{center}
\vspace{-40pt}
\includegraphics[scale=0.55]{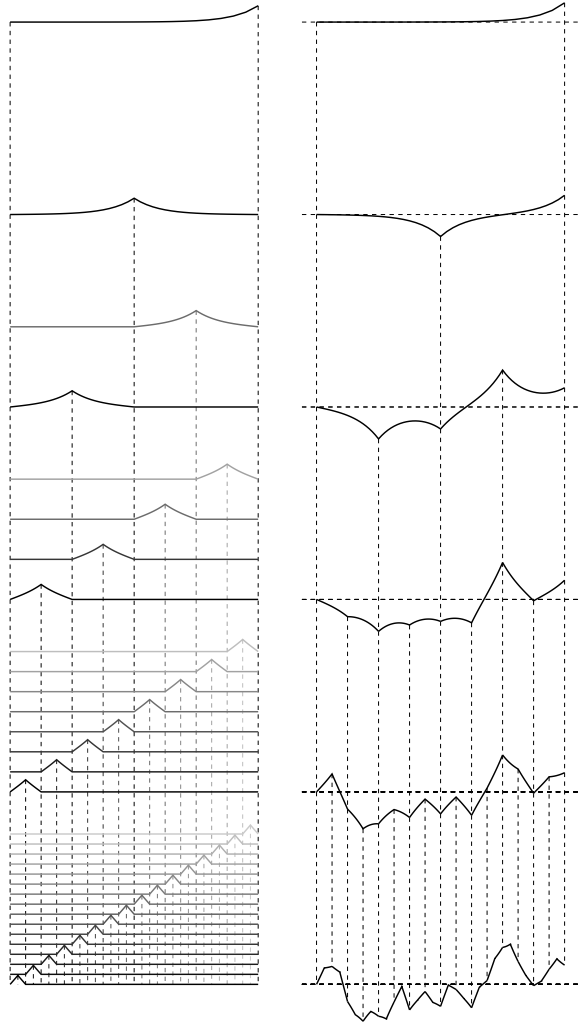}
\end{center}
\caption{In the left column, the elements of the basis $\Phi_{n,k}$ are represented for each rank $n$ with $0 \! \leq \! n \! < \! 6$. In the right column, the conditional Ornstein-Uhlenbeck process $ {\langle U_t \rangle}_{D_n}$ is shown for a given set of realizations on $D_n$.  Once more, note that each element $\Psi_{n,k}$ has a compact support delimited by dyadic numbers in $D_{n} = \lbrace k2^{-n} \, \vert \, 0  \! \leq \! k \! \leq \! 2^{n}\rbrace$ and that all $\Phi_{n',k}$  is zero on $D_{n}$ for $n'>n$.}
\label{fig:OUbasis}
\end{figure}

In view of representing an Ornstein-Uhlenbeck process as a discrete process, the comparison with a Wiener process suggests a candidate basis of decomposition of the form $\Phi_{n,k} \cdot \xi_{n,k}$, the variable $\xi_{n,k}$ following the law $\mathcal{N}(0,1)$.
The deterministic function $\Phi_{n,k}$ is defined with support $S_{n,k} \!= \! \left[ k \! \cdot \!2^{-n\!+\!1}, (k\!+\!1)2^{-n\!+\!1}\right]$ for $n > 0$ with $0 \! \leq \! 2k \! < \! 2^{n}$.
We use expressions \eqref{eq:muOU} and \eqref{eq:sigmaOU} to make explicit the formulation of $\Phi_{n,k}$ in relation \eqref{eq:defPsi} and we obtain
\begin{equation} \label{eq:Phi definition}
\Phi_{n,k}(t)= \left\{ \begin{array}{lllllll}
\displaystyle \sqrt{ \frac{\Gamma}{\alpha}}  \cdot \frac{\sinh{ \left( \alpha  \left( t \! - \! 2k \! \cdot \!2^{-n} \right) \right) }}{\sqrt{\sinh{ \left( \alpha 2^{-n+1} \right) }}}\vspace{5pt}\\
\textrm{if $\left(2k\right)2^{-n} \! \leq \! t \! \leq \! \left(2k\!+\!1\right)2^{-n} $} \: , \\
\\
\displaystyle \sqrt{ \frac{\Gamma}{\alpha}}  \cdot  \frac{\sinh{ \left( \alpha \left( 2\left(k \! + \! 1 \right) 2^{-n} \! - \! t \right) \right) }}{\sqrt{\sinh{ \left( \alpha 2^{-n+1} \right) }}} \vspace{5pt}\\
\textrm{if $\left(2k\!+\!1\right)2^{-n} \! \leq \! t \! \leq \! 2\left(k\!+\!1\right)2^{-n}$} \: ,\\
\\
 0  \quad \textrm{otherwise} \: . \end{array} \right.
\end{equation}
Without any further comment, the element $\Phi_{0,0}$ is defined as
\begin{equation}
\Phi_{0,0}(t) = \sqrt{ \frac{\Gamma}{\alpha}} \cdot \frac{e^{- \frac{\alpha}{2}} \sinh{ \left( \alpha t \right) }}{\sqrt{\sinh{  \alpha  }}} \: ,
\end{equation}
a choice we will explain in the following section.
The first elements $\Phi_{n,k}$ are shown in figure \ref{fig:OUbasis}. As expected, they are only asymptotically scale-invariant
but they exhibit the desirable property of being compactly supported on $S_{n,k}$, the interval between two dyadic points $k \! \cdot \!2^{-n\!+\!1}$ and $(k\!+\!1)2^{-n\!+\!1}$. 
\\
To validate the decomposition of an Ornstein-Uhlenbeck process $U_t$ on the set of  functions $\Phi_{n,k}$, we need to study the convergence of the partial sums
\begin{equation} \label{eq:OU sum 0}
U_{t}^{N}(\omega) = \sum_{n=0}^{N} \sum_{ 0\leq k<2^{n\!-\!1} }\Phi_{n,k}(t) \cdot \xi_{n,k}(\omega) \: .
 \nonumber
\end{equation}
As each function $\Phi_{n,k}$ is dominated by the Haar-derived element $\Psi_{n,k}$, the normal convergence of the series $W_{t}^{N}(\omega)$ in \eqref{eq: Brownian sum} entails the normal convergence of $U_{t}^N(\omega)$ almost surely on the sample space $\Omega$. 
With anticipation of its Ornstein-Uhlenbeck nature, we denote $U_t$ the corresponding limit process.  
The normal convergence causes $U_t$ to be continuous and, being a sum of Gaussian variables, a Gaussian process. 
Therefore, proving that $U_t$ is an Ornstein-Uhlenbeck process just requires us to show that the covariance of $U_t$ satisfies
\begin{equation}\label{eq:covOU}
\langle U_{t} \cdot U_{s } \rangle = \frac{\Gamma}{2 \alpha} \cdot e^{-\alpha \left( t\!+\!s \right)}\left(e^{2 \alpha \left( t \wedge s \right) }-1 \right) \: .
\end{equation}
To establish this relation, we need to evaluate the covariance of $U_t$  as the limit covariance of the partial sums
\begin{equation} \label{eq: covariance OU sum 0}
\langle U_{t}^{N} \cdot U_{s}^{N} \rangle =  
 \sum_{n=0}^{N} \sum_{ 0\leq k<2^{n\!-\!1} }
 \Phi_{n,k}(t) \cdot \Phi_{n,k}(s)  \: .
 \nonumber
\end{equation}
It is possible to simplify the above expression, even though the functions $\Phi_{n,k}$ are not orthogonal. 
For each given $n$, the disjoint supports of $\Phi_{n,k}$ forms a partition of $\left[ 0,1 \right]$ as a collection of segments $S_{n,k}$ of equal length $2^{-n+1}$. 
Considering a real $t$, there is only one sequence of indexes $k_n$ such that $t$ belongs to each support of $S_{n,k_n}$.
The succession of $k_n$ represents $t$ as the intersection of decreasing dyadic segments $ \cap^{\infty}_{n=0} S_{n,k_n}$, which can be explained in terms of the binary representation $t = \sum_{1}^{\infty} a_{i} 2^{-i}  ,   a_{i} \in \lbrace 0,1\rbrace$, if we exclude inappropriate infinite developments.
Bearing in mind the system of indexing for $S_{n,k}$, a simple recurrence argument leads to the expression of $k_n$ corresponding to a given $t$ in its binary representation
\begin{equation}\label{eq:indexKN}
k_n = \frac{1}{2} \cdot \sum_{i=1}^{n-1} a_{i}2^{n-i} \: .
\end{equation}
We are now in a position to write the reduced expression of the partial sums
\begin{equation} \label{eq:OU sum 1}
U_{t}^{N}(\omega) = \sum_{n=0}^{N}\Phi_{n,k_n}(t) \cdot \xi_{n,k_n}(\omega)
\nonumber
\end{equation}
where the terms $\Phi_{k_n,n}(t)$ is made explicit using the previous formulation of $k_n$ in the definition \eqref{eq:Phi definition}
\begin{eqnarray} \label{eq: k_n definition}
\Phi_{n,k_n}(t) =  \left\{ \begin{array}{lll}
\displaystyle \sqrt{ \frac{\Gamma}{\alpha}} \! \cdot \! \frac{\sinh{ \left( \alpha \sum_{n+1}^{\infty} \! a_{i}2^{-i} \right) }}{\sqrt{\sinh{ \left( \alpha 2^{-n+1} \right) }}} & \textrm{if} \; a_{n} \!\! =\!0 \: ,\\
\nonumber\\
\displaystyle  \sqrt{ \frac{\Gamma}{\alpha}} \! \cdot \! \frac{\sinh{ \left( \alpha \sum_{n+1}^{\infty} \! \left( 1\! - \! a_{i} \right) 2^{-i} \right) }}{\sqrt{\sinh{ \left( \alpha 2^{-n+1} \right) }}}  & \textrm{if} \; a_{n} \!\! =\!1 \: .
\end{array} \right.
\end{eqnarray}
Informed by these preliminaries, we shall carry out the calculation of the covariance.
The reduced formulation of partial sums allows us to write
\begin{equation} \label{eq:covarianceOUsum0}
\langle U_{t}^{N} \cdot U_{s}^{N} \rangle =  
\sum_{n=1}^{N} \Phi_{n,k_n}(t) \Phi_{n,l_{n}}(s) + \: \Phi_{0,0}(t) \Phi_{0,0}(s)  \: .
\end{equation}
where the indexes $k_n$ and $l_n$ designate the sequence of functions $\Phi_{k_{n},n}$ and $\Phi_{l_{n},n}$ whose supports contain $t$ and $s$ respectively.
When $t$ and $s$ are distinct, we notice that for $n\!>\!1\!-\! \log_{2}{\vert t\!-\!s \vert}$, the supports $S_{k_n,n}$ and $S_{l_n,n}$ containing $t$ and $s$ respectively are disjoint, so that the cross-products $\Phi_{k_{n},n}(t) \Phi_{l_{n},n}(s)$ cancel out if $n$ is large enough. 
It is then possible to write expression \eqref{eq:covarianceOUsum0} as a finite sum where the terms $\Phi_{k_{n},n}(t)$ and $\Phi_{l_{n},n}(s)$ are specified due to the binary representations $t = \sum_{1}^{\infty} a_{i} 2^{-i}$ and $s = \sum_{1}^{\infty} b_{i} 2^{-i}$. 
We specify that we only consider proper binary representations, that is, the binary representation of dyadic points is chosen in its finite form.
For the sake of simplicity, we assume that $t\!<\!s$. 
Formulated in the binary representation, the order $t\!<\!s$ is equivalent to the existence of a natural $N_0\!>\!0$ such that $a_n\!=\!b_n$ as long as $n\!<\!N_0$ and $a_{N_0}\!<\!b_{N_0}$, that is $a_{N_0}\!=\!0$ and $b_{N_0}\!=\!1$. 
With the preceding remarks, it is clear that  $S_{k_n,n}$ and $S_{l_n,n}$ are disjoint for $n>N_0$ and we can write the covariance of $U_t^N$ for $N>N_0$ in the explicit form
\begin{eqnarray}  \label{eq:covarianceOUsum1}
{\langle U_{t}^{N} \cdot U_{s}^{N} \rangle}_{N>N_0}= 
 \frac{ \Gamma}{2\alpha} \: \Bigg( \sum_{n=1}^{N_0} \frac{2\cdot u_{n}}{\sinh{ \left( \alpha 2^{-n+1} \right) }}  +  e^{-\alpha} \frac{ 2 \sinh{\left( \alpha t \right)} \sinh{\left( \alpha s  \right) } }{\sinh{ \alpha }} \Bigg) \: .
\end{eqnarray}
The variable $u_n$ apparent in  \eqref{eq:covarianceOUsum1} represents  for $n\! < \!N_0$ the numerator of the cross-products $\Phi_{k_{n},n}(t) \Phi_{l_{n},n}(s)$ when the extension $a_n$ and $b_n$ coincide
\begin{eqnarray} \label{eq: u_n definition}
u_{n}= \left\{ \begin{array}{lllll}
\sinh{ \Big( \alpha \sum_{n+1}^{\infty} a_{i}2^{-i} \Big) }\sinh{ \Big( \alpha \sum_{n+1}^{\infty} b_{i}2^{-i} \Big) } \vspace{5pt} \\
 \textrm{if} \quad a_{n} \! =b_{n} \! =0 \: ,\\
\\
\sinh{ \Big( \alpha \sum_{n+1}^{\infty} \left( 1\! - \!  a_{i} \right) 2^{-i} \Big) }\sinh{ \Big( \alpha \sum_{n+1}^{\infty} \left( 1\! - \!  b_{i} \right) 2^{-i} \Big) } \vspace{5pt} \\
\textrm{if} \quad a_{n} \! =b_{n} \! =1 \: ,
\end{array} \right.
\nonumber
\end{eqnarray}
As for the limit case $n=\!N_0$, $u_{N_0}$ expresses the numerator of the cross-product $\Phi_{k_{N_0},N_0}(t) \Phi_{l_{N_0},N_0}(s)$ with $a_{N_0}\!=\!0$ and $b_{N_0}\!=\!1$
\begin{equation}
u_{N_0} = \sinh{ \left( \alpha \sum_{N_0+1}^{\infty} a_{i}2^{-i} \right) } \sinh{ \left( \alpha \sum_{N_0+1}^{\infty}\left(1\!-\! b_{i}\right)2^{-i} \right) } \: .
\nonumber\\
\end{equation}
At that point, the explicit form of the covariance \eqref{eq:covarianceOUsum1} results in a rather complicated combination of hyperbolic functions.
Fortunately enough, we can resort to using remarkable identities to simplify its expression.
The solution actually lies in the consideration of the quantity
\begin{equation} \label{eq: v_n definition}
v_{n} = \sinh{ \left( \alpha \sum_{n}^{\infty} a_{i}2^{-i} \right)} \cdot \sinh{ \left( \alpha \sum_{n}^{\infty} \left( 1\!-\!b_{i} \right) 2^{-i} \right) } \: .
\end{equation}
We show in the supplementary materials that, as long as  $n\!<\!N_0$, $v_n$ verifies the recurrence relation
\begin{equation} \label{eq: v_n recurrence}
v_{n} = 2\cosh \left( \alpha 2^{n} \right) \cdot v_{n+1} + u_{n} \: .
\end{equation}
We can express $u_n$ in terms of $v_n$ and $v_n+1$ to compute the following series by cancellation term by term
\begin{equation}\label{eq:sumDiff}
 \sum_{n=1}^{N_0-1} \frac{ u_n}{\sinh{ \left( \alpha 2^{-n+1} \right) }}  =  \frac{ v_1}{\sinh{ \left( \alpha / 2\right)} } - \frac{ v_{N_0}}{\sinh{ \left( \alpha 2^{-N_0\!+\!1 } \right) }} \: .
\end{equation}
Remembering that $a_{N_0}\!=\!0$ and $b_{N_0}\!=\!1$, we remark that $v_{N_0}\!=\!v_{N_0+1}\!=\!u_{N_0}$ so that the insertion of \eqref{eq:sumDiff} in expression \eqref{eq:covarianceOUsum1} caused the remaining terms in $u_{N_0}$ to cancel out.
 It is then straightforward to write the covariance 
\begin{eqnarray}  \label{eq:covarianceOUsum2}
{\langle U_{t}^{N} \cdot U_{s}^{N} \rangle}_{N>N_0}= 
\frac{ \Gamma}{2\alpha} \: \Bigg(  \frac{2\cdot v_1}{\sinh{ \left( \alpha /2 \right) }}  +  e^{-\alpha} \frac{ 2 \sinh{\left( \alpha t \right)} \sinh{\left( \alpha s  \right) } }{\sinh{ \alpha }} \Bigg) \: ,
\nonumber\\
\end{eqnarray}
We observed that the definition of $v_1$ invokes the full binary representations of $t$ and $s$ so that we have $v_1\! = \! \sinh{ \left( \alpha t \right) } \sinh{ \left( \alpha \left( 1\! - \! s \right)  \right) }$. After a several manipulations, expression \eqref{eq:covarianceOUsum2} finally yields
\begin{equation} \label{eq:covST}
 {\langle U_{t}^{N} \cdot U_{s}^{N} \rangle}_{N>N_0} = 
 \frac{\Gamma}{2 \alpha} \cdot e^{-\alpha \left( t\!+\!s \right)}\left(e^{2 \alpha t }-1 \right) \: ,
 \nonumber
\end{equation}
which is the expected result for the covariance of an Ornstein-Uhlenbeck process \eqref{eq:covOU} given that $t\! = \min \left( t , s \right)$ as $t\!<\!s$.\\
Regarding the calculation of the variance when $t\!=\!s$, the series of cross-products $\Phi_{k_n,n}(t)\Phi_{l_n,n}(s)\!=\!\Phi^{2}_{k_n,n}(t)$ becomes infinite, but fortunately the recurrence relation \eqref{eq: v_n recurrence} is then valid for every $n>0$. 
As the quantity $v_n$ vanishes when $n$ grows to infinity, the cancellation term by term is still effective to compute the series in \eqref{eq:covarianceOUsum1}. It leads to the expected variance expression for an Ornstein-Uhlenbeck process ${}_{ \scriptscriptstyle U \displaystyle} \sigma_{t}^{2} = \frac{\Gamma}{2 \alpha} \! \cdot \! \left( 1 - e^{-\alpha t} \right)$.\\
We finally recap the result for any $t$ and $s$ without assuming any order
\begin{equation} \label{}
\langle U_{t} \cdot U_{s} \rangle = 
\lim_{N \to \infty}  \langle U_{t}^{N} \cdot U_{s}^{N} \rangle = 
 \frac{\Gamma}{2 \alpha} \cdot e^{-\alpha \left( t\!+\!s \right)}\left(e^{2 \alpha \left( t \wedge s \right)  }-1 \right) \: .
 \nonumber
\end{equation}
It proves the discrete description of an Ornstein-Uhlenbeck processes as the normally convergent series of random functions $\Phi_{n,k} \cdot \xi_{n,k}$, where $\Phi_{n,k} $ is a deterministic function
defined in \eqref{eq:Psi definition} and $\xi_{n,k}$ a random variable of normal law $\mathcal{N} \left( 0,1\right)$.


\section{Representation as a bi-infinite sum of random functions}

Whether standing for a Wiener process or an Ornstein-Uhlenbeck process, $X_t$ can be decomposed in a discrete basis of  functions $f_{n,k}$, where $f_{n,k}$ is a generic notation for the deterministic functions  $\Psi_{n,k}$ and $\Phi_{n,k}$.
It suggests to consider the process $X_t$ as a recurrence construction, a view that explains how to chose the first element of the basis $f_{0,0}$.\\
Imagine we want to build a sample path of the continuous process $X_t$ starting with the prior knowledge of its values on the dyadic set $D_N$.
To proceed at the next stage $N\!+\!1$, we need to establish the values of $X_t$ on $D_{N\!+\!1} \! \setminus \! D_{\!N}$, which implies the drawing of as many Gaussian random variables as there are points in this set.
If we consider a given time $t$, there exists a unique $k$ such that $k2^{-N} \! \leq \! t \! < \! (k \! + \! 1)2^{-N}$ and we know that the collection of segments $S_{N\!+\!1, k} \!=\! \lbrack k2^{-N}, (k \! + \! 1)2^{-N} \rbrack$ for $0 \! \leq \! k < \! 2^N$ defines a partition of $\lbrack 0, 1\rbrack$.
We also remark that $t_{N\!+\!1,k} \! = \! (2k\!+\!1)2^{-\left( N\!+\!1\right)}$ is the only point of $D_{N\!+\!1} \! \setminus \! D_{\!N}$ in $S_{N\!+\!1, k}$ and consequently, we note $\zeta_{N\!+\!1,k}$ the Gaussian drawing occurring there. 
According to the results exposed in the second section, we posit 
\begin{equation}
\zeta_{N\!+\!1,k} = {}_{\scriptscriptstyle X \displaystyle}\sigma(t_{N\!+\!1,k}) \cdot \xi_{N\!+\!1,k} + {}_{\scriptscriptstyle X \displaystyle}\mu(t_{N\!+\!1,k}) \: .
\nonumber
\end{equation}
where $\xi_{N\!+\!1,k}$ is of normal law $\mathcal{N}(0,1)$.
Repeating such a construction for $n\!>\!N\!+\!1$ leads us to evaluate the sample path on the whole set of dyadic points $\mathcal{D}$. As $\mathcal{D}$ is dense in $\lbrack 0,1\rbrack$, the complete path is naturally obtained by continued extension. \\
To construct a process rather than a sample path, we need to formulate the above recurrence argument in terms of random functions.
We consider the conditional average ${\langle X_{t} \rangle}_{D_N}$, which is also the most probable outcome of $X_t$ knowing its values on $D_N$.
As $X_t$ is a Markov process, the change in the estimation  ${\langle X_{t} \rangle}_{D_{N\!+\!1}} \!\! - {\langle X_{t} \rangle}_{D_N}$ only depends on the outcome of $\xi_{N\!+\!1,k}$ when restricted on the support $S_{n\!+\!1, k}$.
Due to the simplicity of the situation, it is possible to find an analytical expression of the form $f_{N\!+\!1,k} \cdot \xi_{n,k}$ to describe ${\langle X_{t} \rangle}_{D_{N\!+\!1}} \!\! - {\langle X_{t} \rangle}_{D_N}$ on $S_{N\!+\!1, k}$.
With the help of the so-defined functions $f_{n,k}$, we can introduce the partial sum
\begin{equation}
X^{N}_t =  \sum_{n=0}^{n=N} \sum_{ 0\leq k<2^{n\!-\!1} } f_{n,k}(t)\cdot \xi_{n,k} 
\nonumber
\end{equation}
and relate it to the conditional average ${\langle P_{t} \rangle}_{D_N}$ considered as a random variable.
By definition of $f_{N\!+\!1}$, if $P^{N}_t$ coincides with ${\langle X_{t} \rangle}_{D_N}$,  $P^{N\!+\!1}_t$ equals ${\langle X_{t} \rangle}_{D_{N\!+\!1}}$ at next step. 
Continuing this identification for $n \! > \! N\!+\!1$ shows that the limit process $\lim_{N \to \infty}P^{N}_t$ agrees with the corresponding Ornstein-Uhlembeck process on the index set $\mathcal{D}$. 
Dealing with continuous processes, the density of $\mathcal{D}$ in $\lbrack 0,1\rbrack$ allows us to extent the results on $\lbrack 0,1\rbrack$.
Incidentally, we have an interpretation for the statistical contribution of a component $f_{n,k}(t) \cdot \xi_{n,k}$.
At each step $N$,the function $\sum_{k} f_{N\!+\!1,k}(t) \cdot \xi_{N\!+\!,k}$ represents the uncertainty about $P_t$ that is discarded by the knowledge of its values on $D_{N\!+\!1} \! \setminus \! D_{\!N}$.\\
To validate the recurrence argument, it now remains to verify the initial statement 
\begin{equation}
P^{0}_t =  f_{0,0}(t) \cdot \xi_{0,0} = {\langle P_{t} \rangle}_{D_0} \: 
\nonumber.
\end{equation}
Actually, the need to satisfy this prerequisite enforces how to set the expression of $f_{0,0}$. 
The conditional average ${\langle X_{t} \rangle}_{D_0}$ is a function of the value of $X_t$ on $D_0 \! = \! \lbrace 0, 1\rbrace$. 
By construction the value of $X_t$ in $0$ is assumed to be zero.
We note $Z_{0,0}$ the random function $X_1$ knowing $X_{0}=0$, and we recall that its statistics is given by relations \eqref{eq:probRW}  for a Wiener process and \eqref{eq:probOU} for an Ornstein-Uhlenbeck process respectively. 
With the notation of the second section, we write  ${\langle X_{t} \rangle}_{D_0}$ as a function of $Z_{0,0}$
\begin{equation}
 {\langle X_{t} \rangle}_{D_0} = {}_{ \scriptscriptstyle X \displaystyle} \mu_{0,1}(t, 0, Z_{0,0}) \:  .
 \nonumber
\end{equation} 
It defines a Gaussian random function ${\langle X_{t} \rangle}_{D_0}$ of the form $f_{0,0} \cdot \xi_{0,0}$. 
The dependency of its variance upon $t$ yields the expression of the deterministic part $f_{0,0}$ 
\begin{equation} \label{eq:varfirst}
f_{0,0}(t) = \sqrt{\langle {\langle X_{t} \rangle}_{D_0}^2 \rangle} \: .
\end{equation}
When applied to the Wiener process, the above relation gives the right expression for $\Psi_{0,0}$; when applied to the Ornstein-Ulenbeck process, it gives the already mentioned expression of $\Phi_{0,0}$.\\
Now, we further this recurrence description to show why $X_t$ is naturally represented as a bi-infinite series of random functions.
In that perspective, we extend the definition of the dyadic sets to $D_{N}=\lbrace k2^{-N} \, \vert \, k \in \mathbb{Z} \rbrace$ and we have $D_{-N} = 2^{N} \mathbb{Z} \subset \cdots \subset D_{0} = \mathbb{Z} \subset  \cdots \subset D_{N} = 2^{-N} \mathbb{Z}$.
If we restrain the description of $P_t$ to the index set $\lbrack 0,1\rbrack$, the argument to set the initial step of the recurrence allows us to define a function $f_{-N,0}^{\star}$ so that $f_{-N,0}^{\star} \cdot \xi_{-N,0} = {\langle X_{t} \rangle}_{D_{-N}}$ on $\lbrack 0,1\rbrack$.
The only requirement to adjust \eqref{eq:varfirst} is that $\rbrack 0,1\lbrack$ has no points in $D_{-N}$. 
The usual recurrence construction is then easily adapted to build the process $P_t$ on that segment:
for $n \!>\! -N$, the analytical expressions of $f_{n,k}^{\star}$ are simple extensions of the usual formulas $f_{n,k}$.
In the case of a Wiener process, we make explicit the functions $\Psi_{n,k}^{\star}$ 
\begin{equation} \label{eq:Psi star definition}
\Psi_{n,k}^{\star}(t)= \left\{ \begin{array}{lllllll}
\displaystyle \sqrt{ \frac{\Gamma}{2^{-n+1}}}  \cdot  \left( t  -  2k  \cdot 2^{-n} \right) \vspace{5pt}\\
\textrm{if $\left(2k\right)2^{-n} \! \leq \! t \! \leq \! \left(2k\!+\!1\right)2^{-n} $} \: , \\
\\
\displaystyle \sqrt{ \frac{\Gamma}{2^{-n+1}}}   \cdot  \left( 2\left(k \! + \! 1 \right) 2^{-n}  -  t \right) \vspace{5pt}\\
\textrm{if $\left(2k\!+\!1\right)2^{-n} \! \leq \! t \! \leq \! 2\left(k\!+\!1\right)2^{-n}$} \: ,\\
\\
 0  \quad \textrm{otherwise} \: . \end{array} \nonumber \right.
\end{equation}
for $n \!>\! -N$, and the first element $\Psi_{-N,0}^{\star}$
\begin{equation} \label{eq:first Psi}
\Psi_{-N,0}^{\star}(t)=\sqrt{\frac{\Gamma}{2^{N}}} \! \cdot \! t \: .
\end{equation}
In the case of an Ornstein-Uhlenbeck process, we similarly write the functions $\Phi_{n,k}^{\star}$
\begin{equation} \label{eq:Phi star definition}
\Phi_{n,k}^{\star}(t)= \left\{ \begin{array}{lllllll}
\displaystyle \sqrt{ \frac{\Gamma}{\alpha}} \cdot  \frac{\sinh{ \left( \alpha  \left( t \! - \! 2k \! \cdot \!2^{-n} \right) \right) }}{\sqrt{\sinh{ \left( \alpha 2^{-n+1} \right) }}}\vspace{5pt}\\
\textrm{if $\left(2k\right)2^{-n} \! \leq \! t \! \leq \! \left(2k\!+\!1\right)2^{-n} $} \: , \\
\\
\displaystyle \sqrt{ \frac{\Gamma}{\alpha}}   \cdot  \frac{\sinh{ \left( \alpha \left( 2\left(k \! + \! 1 \right) 2^{-n} \! - \! t \right) \right) }}{\sqrt{\sinh{ \left( \alpha 2^{-n+1} \right) }}} \vspace{5pt}\\
\textrm{if $\left(2k\!+\!1\right)2^{-n} \! \leq \! t \! \leq \! 2\left(k\!+\!1\right)2^{-n}$} \: ,\\
\\
 0  \quad \textrm{otherwise} \: . \end{array}  \nonumber \right.
\end{equation} 
for $n \!>\! -N$, and the first element $\Phi_{-N,0}^{\star}$
\begin{equation} \label{eq:first Phi}
\Phi_{-N,0}^{\star}(t)=\sqrt{\frac{\Gamma}{\alpha}} \cdot  \frac{e^{-\alpha2^{N\!-\!1}} \sinh{ \left( \alpha t \right) }}{ \sqrt{ \sinh{\alpha2^{N}} } } \: .
\end{equation}
As apparent in \eqref{eq:first Psi} and \eqref{eq:first Phi}, the upper bound of $f_{-N,0}$ on $\lbrack 0, 1\rbrack$ is exponentially decreasing toward zero when $N$ goes to infinity.  
The exponential uniform convergence of $f_{-N,0}$ toward zero prescribes to represent the process $X_t$ as a bi-infinite series of random functions $f_{n,k}^{\star} \cdot \xi_{n,k}^{\star}$
\begin{equation} 
X_{t}^{\star N} = \sum_{n=1}^{N} f_{n,k_{n}}^{\star}(t)\xi_{n,k_{n}}^{\star} + \sum_{n=-N}^{0} f_{n,0}^{\star}(t)\xi_{n,0}^{\star} \: . \nonumber
\end{equation}
In the partial sum $X_{t}^{\star N}$, the index $k_{n}$ refers to the unique functions $f_{n,k}^{\star}$ whose support contains $t$ for a fixed $n$. 
In the case of  $n \! \leq \! 0$, this index is constantly set to zero. 
The two series in $X_{t}^{\star N}$ are normally convergent on $\lbrack 0,1 \rbrack$. 
To verify the cogency of the bi-infinite decomposition, it is enough to demonstrate that the covariance of the partial sum 
\begin{equation} \label{eq:covarianceOUdoublesum0}
\langle X_{t}^{\star N} \cdot X_{s}^{\star N} \rangle =  
\sum_{n=1}^{N} f_{n,k_{n}}^{\star}(t) f_{n,l_{n}}^{\star}(s) +  \sum_{n=-N}^{0} f_{n,0}^{\star}(t) f_{n,0}^{\star}(s)\: .
\nonumber
\end{equation}
converges toward the expected covariance of the process $X_t$.
This amounts to show that the right series in the above expressions equals to $f_{0,0}(t) \cdot f_{0,0}(s)$. 
In the case of a Wiener process, a direct calculation leads to
\begin{equation}
\sum_{n=-N}^{0} \Psi_{n,0}^{\star}(t) \Psi_{n,0}^{\star}(s)  = \sum_{n=0}^{N} \frac{\Gamma}{2^{n+1}} \cdot t s = \Gamma \cdot ts \: ,
\nonumber
\end{equation}
which is the exact expression of $\Psi_{n,0}(t) \Psi_{n,0}(s)$. 
In the case of an Ornstein-Uhlenbeck process, a similar calculation yields 
\begin{eqnarray}
\sum_{n=-N}^{0} \Phi_{n,0}^{\star}(t) \Phi_{n,0}^{\star}(s) & = &  \sum_{n=0}^{N} \frac{\Gamma}{\alpha} \cdot \frac{\sinh{\left( \alpha t \right) } \sinh{ \left( \alpha s \right) }}{\sinh{ \left( \alpha 2^{n\!+\!1} \right) }}
\nonumber\\
& = & \frac{\Gamma}{\alpha} \cdot e^{-\alpha}\frac{\sinh{\left( \alpha t \right) } \sinh{ \left( \alpha s \right)}}{\sinh{\alpha}} \: , \nonumber\\
\nonumber
\end{eqnarray}
which consistently equals to $\Phi_{0,0}(t) \Phi_{0,0}(s)$. The previous derivation requires the use of the identity
\begin{equation}
\sum_{n=1}^{\infty} \frac{1}{\sinh{ \left( \alpha 2^{n} \right) }} = \frac{e^{-\alpha}}{\sinh{\alpha}} \nonumber
\end{equation}
proven in the supplementary materials.\\
The representation of a Wiener process as a bi-infinite series of random functions stems from its scale-invariance.
During the construction process, the values of $W_t$ on $D_{N\!+\!1} \! \setminus \! D_{N}$ only depends upon the previous drawings on $D_{N}$. 
Initial drawings on the asymptotic bounds of $\lbrace -\infty, 0, \infty \rbrace = \lim_{N \to -\infty} D_{N}$ only produce vanishing correlations for later stage.
At any finite step $N\!+\!1$ in $\mathbb{Z}$, evaluating $W_t$ on  $D_{N\!+\!1} \! \setminus \! D_{N}$ is a completely scale-invariant operation, which justifies a common analytical expression for the elements $\Psi_{n,k}$.
For nonzero coefficient $\alpha$, the linear component in the Langevin equation precludes any scale-invariance for an Ornstein-Uhlenbeck process $U_t$.
In that respect, it is quite remarkable that we can decompose  $U_t$ in a basis of functions $\Phi_{n,k}$ defined with a unique analytical expression.
If $n \! > \! \log_{2}{\alpha}$ the functions $\Phi_{n,k}$ are well approximated by the functions $\Psi_{n,k}$ showing the typical scale-invariance of a Wiener process; if $n \! < \! \log_{2}{\alpha}$ the functions $\Phi_{n,k}$ are exponential attenuation of the functions $\Psi_{n,k}$, with constant extremal value $\sqrt{\Gamma/2\alpha}$.
As is obvious, the parameter $\alpha$ can be interpreted in terms of a characteristic time.
For time-scales smaller than $2/\alpha$, the components $\Phi_{n,k}$ with intersecting supports add to each others so that the resulting process displays the same correlations as a Wiener process. 
For time-scales larger than $2/\alpha$, the components $\Phi_{n,k}$ add independently  because the exponential attenuation causes all the components but one to be negligible on intersecting supports.
\\
\\
We thank Mariela Sued and Daniel Andor. This work has been supported in part by NIH under grant R01-DC07294.

\end{document}